\documentclass[preprint,12pt]{elsarticle}

\usepackage{amssymb}
\usepackage{amsmath}
\newtheorem{thm}{Theorem}

\newtheorem{lem}[thm]{Lemma}
\newdefinition{rem}{Remark}
\newdefinition{defi}{Definition}
\newproof{pf}{Proof}
\newproof{pot}{Proof of Theorem \ref{thm2}}

\def\Sp{\operatorname{Sp}}
\def\arctanh{\operatorname{arctanh}}

\journal{*****}

\begin{document}
\begin{frontmatter}
\title{Non integrability table for homogeneous potentials invariant by rotation}
\author{Thierry COMBOT\fnref{label2}}
\ead{combot@imcce.fr}
\address{IMCCE, 77 Avenue Denfert Rochereau 75014 PARIS}

\title{Integrability conditions at order $2$ for homogeneous potentials of degree $-1$}
\author{}
\address{}
\begin{abstract}
We prove a meromorphic integrability condition at order $2$ near a homothetic orbit for a meromorphic homogeneous potential of degree $-1$, which extend the Morales Ramis conditions of order $1$. Conversely, we prove that if this criterion is satisfied, then the Galois group of second variational equation is abelian and we compute explicitly the Galois group and the Picard-Vessiot extension.
\end{abstract}
 

\end{frontmatter}

\section{Introduction}

We want to study the Hamiltonian systems of the following form
$$H=T+V$$
with $T=\sum \frac{p_i^2}{2}$ and $V$ a homogeneous function of degree $-1$ in $q_1,\dots,q_n$ and $n\in\mathbb{N}^*$. The corresponding differential equation system is written
$$\dot{q_i}=\frac{\partial}{\partial p_i} H \quad \dot{p_i}=\frac{\partial}{\partial q_i} H \quad\quad (S)$$
For the following and our non-integrability proofs, we will only consider meromorphic potentials $V$.

\begin{defi}
We call $c \in \mathbb{C}^n$ Darboux point a solution of the equation
$$V'(c)=\alpha c$$
with $\alpha \in \mathbb{C}$ called multiplier. Because of homogeneity, we can always choose $\alpha=0,-1$, and we will say that $c$ is non degenerated if $\alpha\neq 0$.
\end{defi}

\begin{defi}
A homothetic orbit associated to a Darboux point $c$ with multiplier $\alpha$ is given by
$$q_i(t)=c_i\phi(t)\quad p_i(t)=c_i \dot{\phi(t)} \quad i=1\dots n$$
with $\phi$ satisfying the following differential equation
$$\frac{1} {2} \dot{\phi(t)}^2=-\frac{\alpha} {\phi(t)}+E$$
with $E$ corresponding the energy orbit.
\end{defi}

In the following, we will call ``norm'' and scalar product the expressions
$$\mid v \mid^2 =\sum\limits_{i=1}^n v_i^2 \qquad <v,w> =\sum\limits_{i=1}^n v_iw_i$$
even for complex $v,w$. We will say moreover that a matrix is orthonormal complex if its columns $X_1,\dots,X_n$ are such that
$$<X_i,X_j>=\sum\limits_{k=1}^n (X_i)_k(X_j)_k=0 \quad\forall i,j\qquad \mid X_i \mid^2=\sum\limits_{k=1}^n (X_i)_k^2=1 \quad\forall i$$

\begin{thm}\label{thmmorales} (Morales, Ramis, Yoshida \cite{1},\cite{4},\cite{5},\cite{6}) Let $V$ be a meromorphic homogeneous potential of degree $-1$ and $c$ a non degenerate Darboux point. If $V$ is meromorphically integrable, then the identity component of the Galois group of the variational equation near the homothetic orbit associated to $c$ is abelian at any order. Moreover, if we fix the multiplier of the Darboux point $c$ to $-1$, then the identity component of the Galois group of the first order variational equation is abelian if and only if
$$\Sp(\nabla^2 V(c))\subset \left\lbrace\frac{1}{2}(k-1)(k+2),\; k\in\mathbb{N} \right\rbrace$$
\end{thm}

The main theorems of this article are the following

\begin{thm}\label{thmmain1}
Let $V$ be a meromorphic homogeneous potential of degree $-1$, $c$ a Darboux point of $V$ with multiplier $-1$. \textbf{We will suppose that $\nabla^2V(c)$, the Hessian of $V$ in $c$, is diagonalizable}. We pose $\lambda_i\;\;i=1\dots n$ its eigenvalues and $X_1,\dots,X_n$ its eigenvectors. If $V$ is meromorphically integrable, then
\begin{enumerate}[i]
\item $\lambda_i=\frac{(p_i-1)(p_i+2)}{2}$ with $p_i\in\mathbb{N}$ (first order integrability condition)
\item $\forall i,j,k=1...n$, $A_{p_i,p_j,p_k}=0 \Rightarrow D^3V(c).(X_i,X_j,X_k)=0$ where $A$ is a $3$ index table with values in $\lbrace 0,1 \rbrace$ invariant by permutation and given by
\begin{itemize}
\item For $i,j,k\in\mathbb{N}^*$, $A_{i,j,k}=1$ if and only if one of the following conditions are satisfied
$$\left\{ \begin{array}{c}
i+j-k\geq 2\\
i-j+k\geq 2\\
-i+j+k\geq 2\\
i+j+k \hbox{ mod } 2=0 \end{array} \right.\hbox{or } \left\{ \begin{array}{c}
-i+j+k\leq -3\\
i+j+k\hbox{ mod } 2=1\end{array} \right.$$
$$\left\{ \begin{array}{c}
i-j+k\leq -3\\
i+j+k\hbox{ mod } 2=1 \end{array} \right. \hbox{or } \left\{ \begin{array}{c}
i+j-k\leq -3\\
i+j+k\hbox{ mod } 2=1 \end{array} \right.$$
\item For $i=0$, $j,k\in\mathbb{N}^*$, $A_{0,j,k}=1$ if and only if $\mid j-k \mid \geq 2$
\item For $i=j=0$, $A_{i,j,k}=1$.
\end{itemize}
\end{enumerate}
\end{thm}

\begin{thm}\label{thmmain2}
Let $V$ be a meromorphic homogeneous potential of degree $-1$, $c$ a Darboux point of $V$ with multiplier $-1$. \textbf{We will suppose that $\nabla^2V(c)$, the Hessian of $V$ in $c$, is diagonalizable}. If $V$ is meromorphically integrable, then the Galois group of the $2$-th order variational equation is always isomorphic to $\mathbb{C}$ and the Picard-Vessiot extension field is
$$\mathbb{C}\left(\phi,\dot{\phi},\ln\left(\frac{1}{2}+\phi\left(1 +\frac{1}{\sqrt{2}}\dot{\phi}\right)\right)\right)$$
except if one (or both) of the two following conditions are satisfied
\begin{itemize}
\item $D^3(V)(c)(v,v,v) \neq 0$ with $\nabla^2V(c) v=-v,\; v\neq 0$
\item $D^3(V)(c)(v,v,w) \neq 0$ with $\nabla^2V(c) v=-v, \; \nabla^2V(c) w=0, \; v,w\neq 0$
\end{itemize}
and in this case the Galois group is $\mathbb{C}^2$ and the Picard-Vessiot field
$$\mathbb{C}\left(\phi,\dot{\phi},\ln\left(\frac{1}{2}+\phi\left(1 +\frac{1}{\sqrt{2}}\dot{\phi}\right)\right),\ln\left( \phi \right)\right)$$
\end{thm}

The first order condition is already known, computed by Yoshida \cite{1} based on classification of hypergeometric functions by Kimura \cite{14}. It has been used many times by \cite{2}, \cite{3}, in particular in the $n$ body problem in the case of homogeneity degree $-1$ in \cite{8}. The Morales Ramis theorem holds for variational equations at any order, and so here we want to study completely the second order, and give an integrability characterization at order $2$. The integrability constraint, when the system of second order variational equations is well written, can be found by computing particular monodromy commutator, like done in \cite{2}. But the true difficulty is not here, it is that we need to study the monodromy for an infinite number of eigenvalues, for all possibles eigenvalues. It comes down to the study of a particular $3$ index sequence, and then to find its zero and non zero entries. Of course it can be easily checked one by one, but this is not enough. By chance, this $3$ index sequence possess an explicit expression, but not easy to find (and to prove). The property behind it is that the monodromy commutator is $D$-finite with respect to the eigenvalue parameters, and so it satisfies a $3$ index linear recurrence with polynomial coefficients. A closed form solution can then be guessed by the gfun package, and then its validity checked. The non nullity can then easily be studied because this closed form expression is a hypergeometric sequence.

\section{Variational equations of order $2$}

\subsection{Reduction of second order variational equations}

For the following, $\phi(t)$ will be a solution of the equation $\frac{1}{2}\dot{\phi}^2=\frac{1}{\phi}+1$. The second order variational equation, if the Hessian matrix is diagonal, can be written
\begin{equation}\label{eq1}
\ddot{X}=\frac{1}{\phi(t)^3} DX +\frac{1}{2}\frac{1}{\phi(t)^4}\left(\begin{array}{c}Y(t)^\intercal T_1 Y(t)\\\cdots\\Y(t)^\intercal T_n Y(t)\\ \end{array}\right)\qquad \ddot{Y}=\frac{1}{\phi(t)^3} DY
\end{equation}
with $D$ diagonal, $T_i \in M_n(\mathbb{C})$. The matrix $D$ is the Hessian matrix of $V$ on the Darboux point $c$, and the matrices $T_i,\; i=1\dots n$ are defined by $T_{i,j,k}=D^3(V)(c).(q_i,q_j,q_k)$. A more detailled construction of second variational equation and higher orders can be found in \cite{17}.

\begin{thm}\label{thm1}
We consider the differential equation~\eqref{eq1}. Suppose that the equation in $Y$ has a virtually abelian Galois group. Let us note $K$ the Picard Vessiot extention of \eqref{eq1}, generated by all solutions of~\eqref{eq1}. The Galois group of~\eqref{eq1} is virtually abelian if and only if all the equations
\begin{equation}\label{eq3}
\ddot{X}_i=\frac{D_{i,i}}{\phi(t)^3} X_i +T_{i,j,k} Y_j(t)Y_k(t) \qquad \ddot{Y_j}=\frac{D_{j,j}}{\phi(t)^3} Y_j\qquad \ddot{Y_k}=\frac{D_{k,k}}{\phi(t)^3} Y_k
\end{equation}
have virtually abelian Galois groups. Moreover, the Galois group depends only on the nullity or non nullity of $T_{i,j,k}$.
\end{thm}

\begin{pf}
Suppose that the Galois group of~\eqref{eq1} is virtually abelian. We write
$$Y(t)=\left(\begin{array}{c} C_1 f_1(t)+C_2 f_2(t)\\\cdots\\C_{2n-1} f_{2n-1}(t)+C_{2n} f_{2n}(t)\\ \end{array}\right)$$
using the fact that $D$ is diagonal. Let us fix all the $C_i=0$ except $C_{2j-1},C_{2j}$ for some fixed $j$. We come down to the equation
$$\ddot{X(t)}=\frac{1}{\phi(t)^3} DX(t) +\frac{1}{2}\frac{1}{\phi(t)^4}\left(\begin{array}{c}T_{1,j,j} Y_j(t)^2\\\cdots\\T_{n,j,j} Y_j(t)^2\\ \end{array}\right)$$
The function $Y_j(t)$ is a fixed function, so it is a non homogeneous linear differential equation. By hypothesis, the Galois group of this equation is virtually abelian. Using the method of variation of the constant, we find that the Galois group depend only on the nullity or non nullity of this equation $T_{i,j,j}$.

We can then substract to $T$ all its terms of the form $T_{i,j,j}$, and still the Galois group of~\eqref{eq1} will be virtually abelian. Now we do the same procedure for the terms in $Y_k(t)Y_j(t)$. We fix all the $C=0$ except for $C_{2j-1},C_{2j}$ and $C_{2k-1},C_{2k}$. As we have $T_{i,j,j}=0$ and $T_{i,k,k}=0$ we get the equation
$$\ddot{X}=\frac{1}{\phi(t)^3} DX +\frac{1}{2}\frac{1}{\phi(t)^4}\left(\begin{array}{c}T_{1,j,k} Y_j(t)Y_k(t)\\\cdots\\T_{n,j,k} Y_j(t)Y_k(t)\\ \end{array}\right)$$
The functions $Y_j(t),Y_k(t)$ are fixed, so this is a non homogeneous linear differential equation. With hypothesis, the Galois group of this equation is virtually abelian.  Using the method of variation of the constant, we find that the Galois group depend only on the nullity or non nullity of this equation $T_{i,j,j}$.
Conversely, the equations~\eqref{eq3} are non homogeneous linear differential equations, so we can sum the solutions of the equations
$$\ddot{X}_i=\frac{D_{i,i}}{\phi(t)^3} X_i +T_{i,j,k} Y_j(t)Y_k(t)$$
to produce all the solutions of equation~\eqref{eq1}.
\end{pf}

The following lemma will have a primary importance in computation of monodromy. In fact, it will be necessary only to compute some sort of residue

\begin{thm}\label{thm2}
We consider $F\in\mathbb{C}(z_1)\left[z_2 \right]$ and
$$f(t)=F\left(t,\arctanh\left(\frac{1}{t}\right)\right)$$
We consider the differential field and the Galois group
$$K=\mathbb{C}\left(t,\arctanh\left(\frac{1}{t}\right),\int f dt \right) \qquad G=\sigma(K,\mathbb{C}(t))$$
If $G$ is abelian, then
\begin{equation}\label{eq2}
\frac{\partial}{\partial \alpha} \underset{t=\infty}{\hbox{Res}}\; F\left(t,\arctanh\left(\frac{1}{t}\right)+\alpha\right) =0\quad \forall \alpha\in\mathbb{C}
\end{equation}
\end{thm}

\begin{pf}
First we recall that if the Galois group $G$ is abelian, then so is the monodromy group, because the monodromy group is always included inside the Galois group. We consider two paths, the eight path $\sigma_1$ around the singularities $-1,1$ and the path around $\sigma_2$. At infinity, $F\left(t,\arctanh\left(\frac{1}{t}\right)+\alpha\right)$ has the following series expansion
$$\int F\left(t,\arctanh\left(\frac{1}{t}\right)+\alpha\right) dt=\sum\limits_{n=n_0}^{\infty} a_n(\alpha)t^n+r(\alpha)ln\; t$$
because the function $\arctanh\left(\frac{1}{t}\right)$ is smooth at infinity. We consider now the commutator
$$\sigma=\sigma_2^{-1}\sigma_1^{-\frac{\beta}{2i\pi}}\sigma_2\sigma_1^{\frac{\beta}{2i\pi}}$$
We have that $\sigma_1^{\frac{\beta}{2i\pi}}(f)=F\left(t,\arctanh\left(\frac{1}{t}\right)+\beta\right)$ and $\sigma_2(ln\; t)=ln\; t +2i\pi$. We conclude that
$$\sigma(f)=f+r(\beta)-r(0)$$
This $r(\alpha)$ correspond to the residue of $F\left(t,\arctanh\left(\frac{1}{t}\right)+\alpha\right)$ at infinity. If the monodromy is abelian, then $\sigma$ should act trivially on $f$. This is the case if and only if $r(\beta)-r(0)\;\;\forall \beta\in \mathbb{Z}$. The function $r$ is polynomial in $\beta$, then $r(\beta)-r(0),\;\;\forall\beta\in\mathbb{C}$. This gives the formula \eqref{eq2}
\end{pf}

\begin{figure}
\begin{center}
\includegraphics[width=5.9cm]{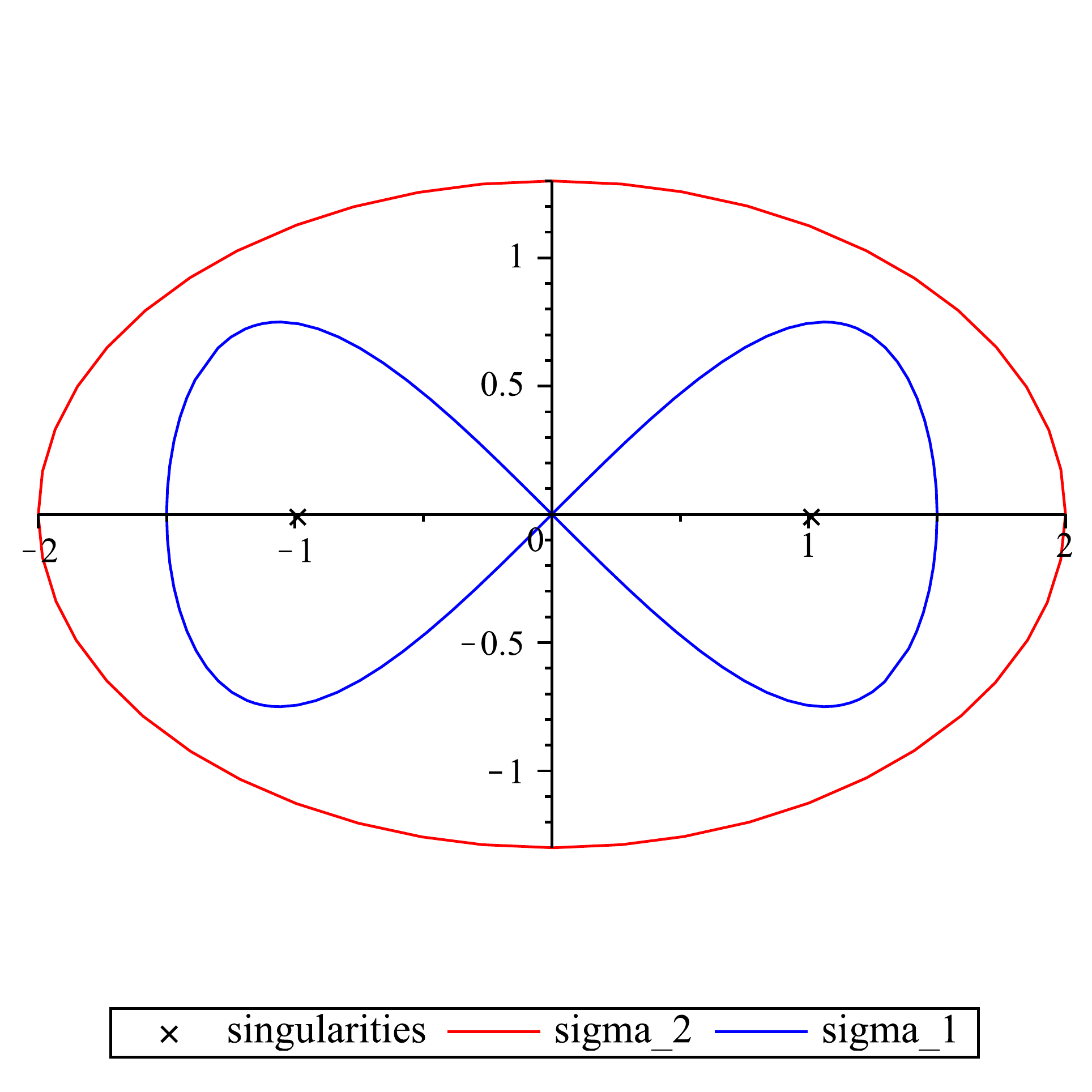}
\includegraphics[width=5.9cm]{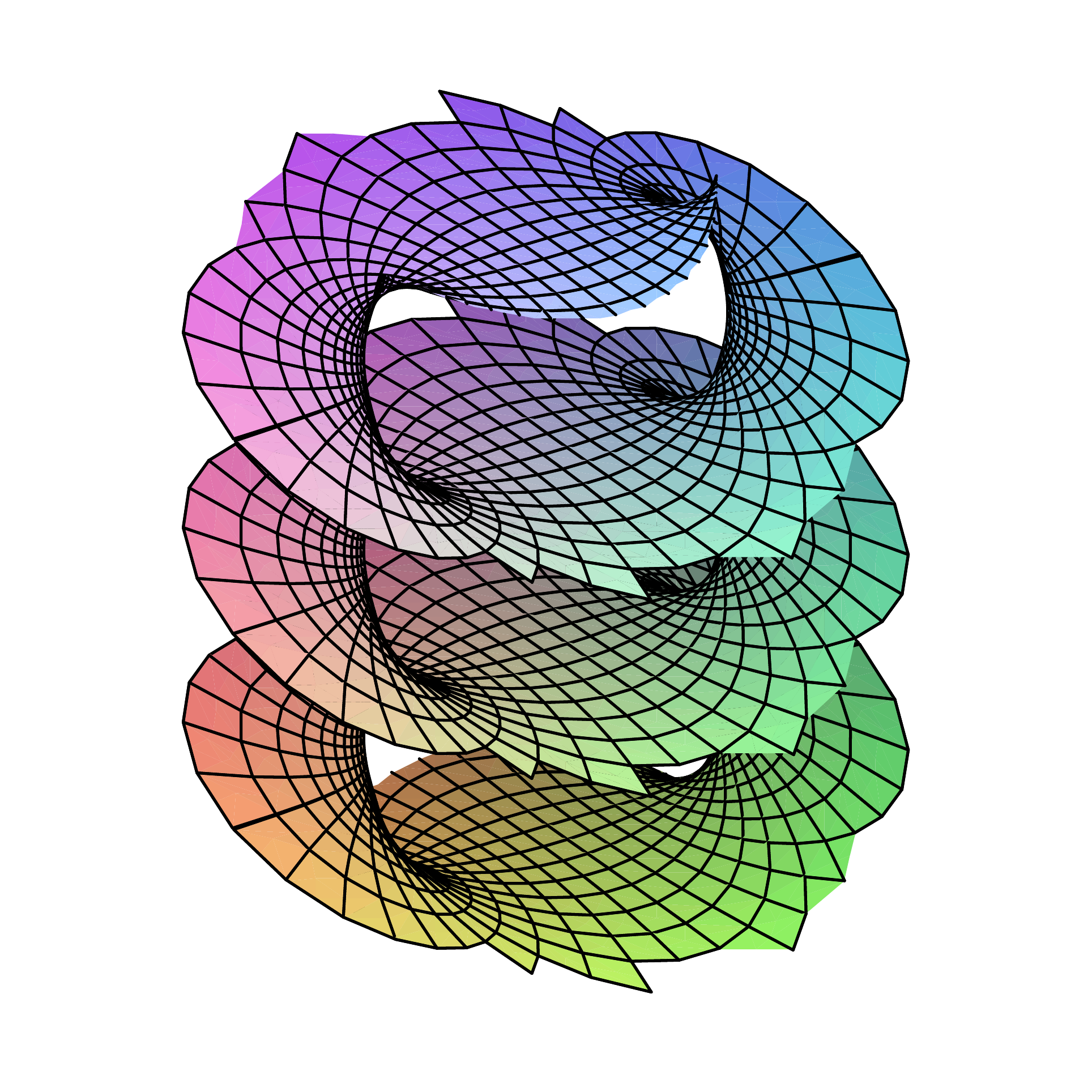}
\end{center}
\caption{Paths corresponding to monodromy elements $\sigma_1,\sigma_2$, and the Riemann surface associated to the $Q_j$. The difference between two sheeves is $i\pi \epsilon_j P_j$. We see that $\sigma_2$, corresponding to monodromy around infinity, acts trivially on $Q_j$.}
\end{figure}

\subsection{Computation of terms of order $2$}

Now we prove that we can always come down to study the equation~\eqref{eq3} to prove Theorems~\ref{thmmain1},~\ref{thmmain2}.

\begin{thm}\label{thm3}
Let $A\in M_n(\mathbb{C})$ be a symmetric matrix, which mean $A_{i,j}=A_{j,i}$. We suppose that $A$ is diagonalizable. Then $A$ is diagonalizable in an orthonormal basis, it means that it exists $X_1,\dots,X_n$ such that
$$<X_i,X_j>=\sum\limits_{k=1}^n (X_i)_k(X_j)_k=0 \quad\forall i,j \qquad \mid X_i \mid^2=\sum\limits_{k=1}^n (X_i)_k^2=1 \quad\forall i$$
\end{thm}

This small result is in fact very important for reduction of problems dealing with homogeneous potentials. This is because the integrability status of a potential $V$ is not changed after a rotation (even a complex one) nor dilatation. Such property allows a great simplification in classification like in \cite{2}, \cite{3}. The simplification at order $2$ is even more important, as a such simple criterion of Theorem~\ref{thmmain1} is not possible if the Hessian matrix is not diagonal. As we will see, the only hypothesis of diagonalizability is very weak, because in particular the system is rarely integrable at order $1$ if the Hessian matrix is not diagonalizable (the conditions on the spectrum are still necessary, but not sufficient).

\begin{pf}
Let $A\in M_n(\mathbb{C})$ be a symmetric diagonalizable matrix, $X_1,\dots,X_n$ and its eigenvectors. Suppose first that $\mid X_k \mid^2\neq 0 \;\;\forall k$. We have then
$$<X_k,Av>=\sum\limits_i (X_k)_i \sum\limits_j a_{i,j} v_j=\sum\limits_j v_j \sum\limits_i a_{i,j}(X_k)_i=$$
$$<A^\intercal X_k,v>=<A X_k,v>=\lambda_k <X_k,v>$$
Then the orthogonal subspace $E_k$ to $X_k$ is stable by $A$. Moreover, $X_k$ is not in $E_k$ then $E_k$ is supplementary to $\mathbb{C}.X_k$, then $E_k$ containts all the other eigenvectors. Then
$$<X_i,X_j>=\sum\limits_{k=1}^n (X_i)_k(X_j)_k=0 \quad\forall i,j$$
Now suppose it exists $X_k$ such that $\mid X_k \mid^2=0$. Then $X_k$ is orthogonal to all eigenvectors with non zero norm (using the proof before). So if it is the only eigenvector with zero norm, then $E_k$ containts all the other eigenvectors, but also containts $X_k$ because $<X_k,X_k>=0$. But $E_k$ is of codimension $1$, then it cannot containt all the eigenvectors (because $X_1,\dots,X_n$ form a basis). Then it exists another vector, let say $X_l$, with zero norm, and moreover $<X_k,X_l>\neq 0$ because otherwise the previous reasoning still holds.\\
Let us now write the matrix $P\in M_n(\mathbb{C})$ formed by the eigenvectors $X_1,\dots,$ $X_n$. We have then
$$A=PDP^{-1}=A^\intercal=(P^{-1})^\intercal D P^\intercal$$
where $D$ is a diagonal matrix. Then
$$P^\intercal P D=D P^\intercal P$$
Then $P^\intercal P$ commute with $D$. We know that $\left[P^\intercal P \right]_{k,l} \neq 0$, then $D_{k,k}=D_{l,l}$ (because $P^\intercal P$ should commute with $D$). Then $X_k$ and $X_l$ have the same eigenvalue. So, the fact that these vectors have zero norm only correspond to a bad ``choice'' of eigenvector. It suffice to take $X_k+X_l$ and $X_k-X_l$ instead of $X_k,X_l$.\\
Now that $\mid X_k \mid^2\neq 0 \;\;\forall k$, we just need to multiply $X_k$ with some constant such that $\mid X_k \mid^2=1 \;\;\forall k$. This end the proof of the theorem.
\end{pf}

Thanks to this theorem, we know that if the Hessian matrix is diagonalizable, we can always make an orthonormal variable change in the potential $V$ (this does not affect at all integrability) such that the Hessian matrix becomes diagonal, and thus apply Theorem~\ref{thm1}. Let us check now that the condition of Theorem~\ref{thmmain1} is well defined and equivalent to this.

\begin{thm}\label{thm4}
Let $V$ be a meromorphic homogeneous potential of degree $-1$ and $c\in\mathbb{C}^n$ a Darboux point with multiplier $-1$. We suppose $\nabla^2V(c)$ diagonalizable and we note its eigenvectors $X_1,\dots,X_n$. Then the integrability constraints of Theorem~\ref{thmmain1} do not depend on the choice of $X_1,\dots,X_n$. Moreover, if we make an orthonormal choice for the $X_1,\dots,X_n$, then coefficients $T_{i,j,k}$ in~\eqref{eq1} will correspond to $D^3(V).(X_i,X_j,X_k)$.
\end{thm}

\begin{pf} 
First, if we make the orthonormal choice (which is always possible thanks to Theorem~\ref{thm3}), we note $P$ the associated orthonormal matrix. We have then that the potential $W(q)=V(P q)$ has a Darboux point in $P^{-1} c$ and the corresponding Hessian matrix is diagonal. Moreover, we have that
$$D^3(W).(q_i,q_j,q_k)=D^3(V).(Pq_i,Pq_j,Pq_k)=D^3(V).(X_i,X_j,X_k) $$
Now, we verify that the criterion is well defined. We first consider the case where all eigenvalues of $\nabla^2 V(c)$ are distinct. Then, up to multiplication by a non zero constant, there is a unique choice of eigenvectors $X_1,\dots,X_n$. So this does not change nullity or non nullity of $D^3(V).(X_i,X_j,X_k)$. Now, if there are multiple eigenvalues. There are an infinite number of choices for eigenvectors $X$. We fix one. Suppose that $X_1,X_2$ have the same eigenvalue. Since the nullity condition is associated only to the corresponding eigenvalues, if there is a nullity condition for some third order derivative involving $X_1$, it will be the same for $X_2$. Suppose there is a condition
$$D^3(V).(X_1,X_j,X_k)=0 \quad D^3(V).(X_2,X_j,X_k)=0$$
Then, for the vector $\alpha X_1+\beta X_2$, we will also have $D^3(V).(\alpha X_1+\beta X_2,X_j,X_k)$ $=0$ by expanding it. Remark that even if $j=1,k=1$, it will still work. All basis changes can be written as successive such linear combinations. So the constraints of Theorem~\ref{thmmain1} do not depend on the choice of $X_1,\dots,X_n$.
\end{pf}

This theorem is in principle not very important, because we can reduce the potential directly to have a diagonal Hessian matrix, but it is useful in pratice. This is because it is not so easy to compute the orthonormal basis change as just to compute the spectrum. In particular when there are parameters, the normalization of eigenvectors becomes uneasy because this produces singularities. We just say here that we can avoid that, and do not care about normalization (for some specific values of the parameters, $X_1,\dots,X_n$ will no more be a basis, but the criterion is in this case still a necessary one).

\section{Non integrability of second order variational equations}

\subsection{A first approach}
Integrability condition at order $2$ for equation~\eqref{eq1}
\begin{thm}\label{thm5}
We consider an equation of the form
$$\ddot{X(t)}=\frac{1}{\phi(t)^3} DX(t) +\frac{1}{2}\frac{1}{\phi(t)^4}\left(\begin{array}{c}Y(t)^\intercal T_1 Y(t)\\\cdots\\Y(t)^\intercal T_n Y(t)\\ \end{array}\right)\qquad \ddot{Y}=\frac{1}{\phi(t)^3} DY $$
with $D$ diagonal, $T_i \in M_n(\mathbb{C})$. The equation~\eqref{eq1} has a virtually abelian Galois group if and only if
\begin{itemize}
\item $D_{i,i}=\frac{(p_i-1)(p_i+2)}{2}$ with $p_i\in\mathbb{N}$ (integrability condition of order $1$)
\item $\forall i,j,k=1...n$, $A_{p_i,p_j,p_k}=0 \Rightarrow T_{i,j,k}=0$ where $A$ is a three index table with values in $\lbrace 0,1 \rbrace$, invariant by permutation and whose first values are given by
\end{itemize}
\end{thm}

\vspace*{-0.1cm}
\begin{tabular}{|c|c|c|c|c|c|c|c|c|c|c|c|c|c|c|c|c|c|}
\hline
$A_{0,i,j}$&0&1&2&3&4&5&6&7&$A_{1,i,j}$&0&1&2&3&4&5&6&7\\\hline
0&1&1&1&1&1&1&1&1&0&1&0&0&1&1&1&1&1\\\hline
1&1&0&0&1&1&1&1&1&1&0&0&0&0&0&1&0&1 \\\hline
2&1&0&0&0&1&1&1&1&2&0&0&0&0&0&0&1&0 \\\hline
3&1&1&0&0&0&1&1&1&3&1&0&0&0&0&0&0&1 \\\hline
4&1&1&1&0&0&0&1&1&4&1&0&0&0&0&0&0&0 \\\hline
5&1&1&1&1&0&0&0&1&5&1&1&0&0&0&0&0&0 \\\hline
6&1&1&1&1&1&0&0&0&6&1&0&1&0&0&0&0&0 \\\hline
7&1&1&1&1&1&1&0&0&7&1&1&0&1&0&0&0&0 \\\hline
$A_{2,i,j}$&0&1&2&3&4&5&6&7&$A_{3,i,j}$&0&1&2&3&4&5&6&7\\\hline
0&1&0&0&0&1&1&1&1&0&1&1&0&0&0&1&1&1\\\hline
1&0&0&0&0&0&0&1&0&1&1&0&0&0&0&0&0&1 \\\hline
2&0&0&1&0&0&0&0&1&2&0&0&0&1&0&0&0&0 \\\hline
3&0&0&0&1&0&0&0&0&3&0&0&1&0&1&0&0&0 \\\hline
4&1&0&0&0&1&0&0&0&4&0&0&0&1&0&1&0&0 \\\hline
5&1&0&0&0&0&1&0&0&5&1&0&0&0&1&0&1&0 \\\hline
6&1&1&0&0&0&0&1&0&6&1&0&0&0&0&1&0&1 \\\hline
7&1&0&1&0&0&0&0&1&7&1&1&0&0&0&0&1&0 \\\hline
$A_{4,i,j}$&0&1&2&3&4&5&6&7&$A_{5,i,j}$&0&1&2&3&4&5&6&7\\\hline
0&1&1&1&0&0&0&1&1&0&1&1&1&1&0&0&0&1\\\hline
1&1&0&0&0&0&0&0&0&1&1&1&0&0&0&0&0&0 \\\hline
2&1&0&0&0&1&0&0&0&2&1&0&0&0&0&1&0&0 \\\hline
3&0&0&0&1&0&1&0&0&3&1&0&0&0&1&0&1&0 \\\hline
4&0&0&1&0&1&0&1&0&4&0&0&0&1&0&1&0&1 \\\hline
5&0&0&0&1&0&1&0&1&5&0&0&1&0&1&0&1&0 \\\hline
6&1&0&0&0&1&0&1&0&6&0&0&0&1&0&1&0&1 \\\hline
7&1&0&0&0&0&1&0&1&7&1&0&0&0&1&0&1&0 \\\hline
$A_{6,i,j}$&0&1&2&3&4&5&6&7&$A_{7,i,j}$&0&1&2&3&4&5&6&7\\\hline
0&1&1&1&1&1&0&0&0&0&1&1&1&1&1&1&0&0\\\hline
1&1&0&1&0&0&0&0&0&1&1&1&0&1&0&0&0&0 \\\hline
2&1&1&0&0&0&0&1&0&2&1&0&1&0&0&0&0&1 \\\hline
3&1&0&0&0&0&1&0&1&3&1&1&0&0&0&0&1&0 \\\hline
4&1&0&0&0&1&0&1&0&4&1&0&0&0&0&1&0&1 \\\hline
5&0&0&0&1&0&1&0&1&5&1&0&0&0&1&0&1&0 \\\hline
6&0&0&1&0&1&0&1&0&6&0&0&0&1&0&1&0&1 \\\hline
7&0&0&0&1&0&1&0&1&7&0&0&1&0&1&0&1&0 \\\hline
\end{tabular}


One direct application of this theorem is the study of problems in celestial mechanics, like in \cite{8},\cite{9}. It is often unnecessary to know the table $A$ for arbitrary high eigenvalues, except in some cases like the open problem at the end of \cite{3}, which is because of this much more difficult.

\begin{pf}
Using Theorem~\ref{thm1}, the study of the Galois group of~\eqref{eq1} comes down to the study of \eqref{eq3}. Let us look first at the expressions of functions $Y$. We can fix $E=1$ because we can fix $E$ to any value using homogeneity (except for $E=0$). The variable change $\phi(t)\longrightarrow t$ gives the equation
$$2t^2(1+t)\ddot{Y}_i-t\dot{Y}_i=\frac{(i-1)(i+2)}{2} Y_i$$
We then make the variable change
$$\sqrt{\frac{1+t}{t}} \longrightarrow t$$
The equations satisfied by $Y_i(t)$ becomes
$$\frac{1}{2}\left( {t}^{2}-1 \right)\ddot{y} +2\,t\dot{y} -\frac{\left( {{i}}-1 \right)\left( {{i}}+2 \right)}{2} y=0$$
A basis of solutions is given by $(P_i,Q_i)$ where $P_i$ are polynomials and the functions $Q_i$ can be written
$$Q_i(t)=P_i(t)\int \frac{1}{(t^2-1)^2P_i(t)^2} dt$$
The functions $Q$ are multivalued except for $i=0$, which is particular, because the Galois group is $Id$ instead of $\mathbb{C}$ and then all solutions are algebraic. We get for $X(t)$ the following solution
$$X(t)=C_1P_i(t)+C_2Q_i(t)+$$ $$\int Y_j(t)Y_k(t)P_i(t)(t^2-1)^2 dtQ_i(t)-\int Y_j(t)Y_k(t)Q_i(t)(t^2-1)^2 dtP_i(t)$$
So we need to study the monodromy of
$$\int Y_j(t)Y_k(t)P_i(t)(t^2-1)^2 dtQ_i(t)-\int Y_j(t)Y_k(t)Q_i(t)(t^2-1)^2 dtP_i(t)$$
Let us try to apply Theorem~\ref{thm2}. This theorem do not apply directly because there could be compensations between the two integrals. But we can rewrite it
$$\int Q_j(t)Q_k(t)P_i(t)(t^2-1)^2 dt P_i(t)\int \frac{1}{(t^2-1)^2P_i(t)^2} dt-$$
$$\int Q_j(t)Q_k(t)P_i(t)\int \frac{1}{(t^2-1)^2P_i(t)^2} dt (t^2-1)^2 dtP_i(t) =$$
$$\int \int Q_j(t)Q_k(t)P_i(t) dt\frac{1}{(t^2-1)^2P_i(t)^2} (t^2-1)^2 dtP_i(t)$$
Then
$$\int Q_j(t)Q_k(t)P_i(t) dt$$
is in the Picard Vessiot field of~\eqref{eq3} (because $P_i$ is a polynomial and that the Picard Vessiot field is stable by derivation). We also have that $Q_i$ is in the Picard Vessiot field, and then by substraction,
$$\int Q_j(t)Q_k(t)Q_i(t)(t^2-1)^2 dt$$
is in the Picard Vessiot field of~\eqref{eq3}. We can now apply Theorem~\ref{thm2} to this integral, and so it is only needed to study the residue
$$S=\underset{t=\infty}{Res}\;  (t^2-1)^2(Q_i(t)+\epsilon_i \alpha P_i)(Q_j(t)+\epsilon_j \alpha P_j)(Q_k(t)+\epsilon_k \alpha P_k) dt$$

The polynomials $P_i$ can be generated by the formula
$$P_i(t)=\frac{1}{t^2-1}\frac{\partial^{i-1}}{\partial t^{i-1}} (t^2-1)^i$$
(which gives a normalization for the dominant coefficient of $P_i$ that we will choose for now) and the functions $Q_i$ can be written
$$Q_i(t)=\epsilon_i P_i(t) \arctanh\left(\frac{1}{t}\right)+\frac{W_i(t)}{t^2-1}$$
with $W_i$ polynomials, $\epsilon_i$ a real sequence.

The integrability at order $2$ only require that the identity component of the Galois group be abelian. This does not say a priori anything about the whole Galois group, except if it is connected. This is the case here. The Galois group at order $1$ is $\mathbb{C}$ (or $Id$). Then at order $2$, the Picard Vessiot field will be of the form
$$K=\mathbb{C}\left(t,\arctanh\left( \frac{1}{t} \right),\int f(t) dt \right) \hbox{ with } f(t)\in \mathbb{C}\left(t, \arctanh\left( \frac{1}{t} \right)\right)$$
We just add some integral in $K$. Then the Galois group $\sigma(K,\mathbb{C}(t))$ is still connected.

So, we write in the table $A_{i,j,k}$, $1$ is the constaint of Theorem~\ref{thm2} is satisfied, $0$ otherwise. This criterion is a priori only a necessary criterion, not sufficient. So we then try to compute the integral using integration by parts
$$\int Q_j(t)Q_k(t)P_i(t)(t^2-1)^2 dtQ_i(t)-\int Q_j(t)Q_k(t)Q_i(t)(t^2-1)^2 dtP_i(t)$$
Using the expressions of functions $Q$, we need to integrate functions in $\mathbb{C}\left[t,\arctanh\left( \frac{1}{t} \right) \right]$. We make successive integration by part, by deriving the term in $\arctanh\left( \frac{1}{t} \right)$ of highest degree, and finally we arrive on an integral of a function in $\mathbb{C}(t)$. This procedure could fail, but it works every time for $A_{i,j,k}=1$ (for $A_{i,j,k}=0$, the procedure fails because terms in $\ln(t^2-1)$ appear). Using Theorem~\ref{thm1}, we know that it is necessary and sufficient that all equations
$$\ddot{X}=\frac{1}{\phi(t)^3} d_i X +\frac{1}{2}\frac{1}{\phi(t)^4}Y_j(t)Y_k(t)$$
have a virtually abelian Galois group for all non zero $T_{i,j,k}$ for the virtual abelianity of the Galois group of~\eqref{eq1}.
\end{pf}

\begin{thm}\label{thm6}
The table $A$ of Theorem~\ref{thm5} has the following values
\begin{itemize}
\item For $i,j,k\in\mathbb{N}^*$, $A_{i,j,k}=1$ if and only if one of the following conditions are satisfied
$$\left\{ \begin{array}{c}
i+j-k\geq 2\\
i-j+k\geq 2\\
-i+j+k\geq 2\\
i+j+k \hbox{ mod } 2=0 \end{array} \right.\hbox{or } \left\{ \begin{array}{c}
-i+j+k\leq -3\\
i+j+k\hbox{ mod } 2=1\end{array} \right.$$
$$\left\{ \begin{array}{c}
i-j+k\leq -3\\
i+j+k\hbox{ mod } 2=1 \end{array} \right. \hbox{or } \left\{ \begin{array}{c}
i+j-k\leq -3\\
i+j+k\hbox{ mod } 2=1 \end{array} \right.$$
\item For $i=0$, $j,k\in\mathbb{N}^*$, $A_{0,j,k}=1$ if and only if $\mid j-k \mid \geq 2$
\item For $i=j=0$, $A_{i,j,k}=1$.
\end{itemize}
Moreover, the table $A$ is invariant by permutation of the index $i,j,k$.
\end{thm}

This table is the direct generalization of the integrability table of \cite{1} at order $2$ for degree $-1$. A similar process could be done for other homogeneity degrees, but in fact the degree $-1$ is much more simple for three reasons
\begin{itemize}
\item There is only one family in the Morales Ramis table for degree $-1$, and generically, there are two (and the complexity increase with the power three of the number of families).
\item Some homogeneity degrees have very particular families, associated the groups $A_4,S_4,A_5$. This produce very very complicated computations.
\item By studying only one homogeneity degree, we have one parameter less in the $D$-finite computations. This is important because computational cost usually increase exponentially with the number of parameters (at least).
\end{itemize}

\subsection{Study of the solutions}

\begin{pf}
Using the last theorem, we already know that we just have to study the residue
$$\underset{t=\infty}{Res}\; (t^2-1)^2(Q_j(t)+\epsilon_j \alpha P_j)(Q_k(t)+\epsilon_k \alpha P_k)(Q_k(t)+\epsilon_k \alpha P_k)$$
and a necessary condition for integrability is that this residue should be independent of $\alpha$. We will call the fact that the coefficient in $\alpha^2$ should be zero ``the constraint in $\alpha^2$'' and respectively ``the constraint in $\alpha$'' for the term in $\alpha$. We see also that a priori, the residue is a polynomial of degree $3$ in $\alpha$ but we have that
$$\hbox{coeff}\left( \underset{t=\infty}{Res}\;  (Q_i(t)+\epsilon_i \alpha P_i)(Q_j(t)+\epsilon_j \alpha P_j)(Q_k(t)+\epsilon_k \alpha P_k)(t^2-1)^2,\alpha^3 \right)=$$
$$\underset{t=\infty}{Res}\;  P_iP_jP_k\epsilon_i\epsilon_j\epsilon_k(t^2-1)^2=0$$
because it is a polynomial in $t$. So we have only two constraints for integrability. Let us begin by checking that all functions $Q_i$ are multivalued for $i\neq 0$. We only need to prove

\begin{lem}
We have $\epsilon_i \neq 0\;\forall i\in\mathbb{N}^*$
\end{lem}

\begin{pf}
Looking at first values of $\epsilon_i$, we already can guess the expression of $\epsilon$
$$\epsilon_i=\frac{4^{-i}i(i+1)}{i!^2}\;i\in\mathbb{N}$$
We now need to prove it. The sequence $\epsilon_i$ can be computed thanks to the formula
$$\epsilon_i=\int\limits_C \frac{1}{(t^2-1)^2P_i(t)^2} dt$$
with $C$ a circle around $-1,1$ in the direct way (because $\epsilon_i$ is the term in front of $\arctanh\left(\frac{1}{t}\right)$ which grows by $1$ along $C$). Using the symmetry $t\longrightarrow -t$, we only need to compute the residue in $1$ for example. We have then
$$\epsilon_i=2\left(\frac{\partial}{\partial t}\left(\frac{1}{(t+1)^2 P_i(t)^2} \right) \right)\mid_{t=1}$$
knowing that $1$ is never a root of $P_i$. So we just need to compute the sequences
$$P_i(1),\;\left(\frac{\partial}{\partial t}P_i(t) \right)\mid_{t=1}$$
with the recurrence formula $P_i$
$$(4n^3+12n^2+8n)P_n+(-4tn^2-14tn-12t)P_{n+1}+(n+3)P_{n+2}$$
We have then
$$P_i(1)=2^i(i+1)!\quad \left(\frac{\partial}{\partial t}P_i(t) \right)\mid_{t=1}=\frac{1}{4} 2^i i (i+3) (i+1)!$$
and we get $\epsilon_i$.
\end{pf}

For the following, we will use the system 
\begin{align*}
\left\lbrace(4n^3+12n^2+8n) f_n(t) - (4tn^2+14tn+12t) f_{n+1}(t) + (n+3) f_{n+2}(t),\right.\\
\left.(t^2-1)f_n''(t) +4tf_n'(t)-(n-1)(n+2)f_n(t)\right\rbrace
\end{align*}
which vanish for $f_n(t)=P_n$ and $f_n(t)=\epsilon_n^{-1} Q_n$. The system
$$\left\{-4tf_n(t)+(t^2-1)f_n'(t),-f_n(t)+f_{n+1}(t)\right\} $$
vanish for $f_n(t)=(t^2-1)^2$. We will use these systems for the package Mgfun to compute recurrences for our residues.

\subsection{General case}

\begin{pf}
First part. We prove that the Galois group is not virtually abelian if the variational equation contains a term corresponding to index such that $A_{i,j,k}=0$. We will begin with the non zero index case (this is because the index $0$ is very special, in particular the function $Q_0$ is not multivalued). We now need to compute the residues of Theorem~\ref{thm2} for all index and prove they are non zero for $A_{i,j,k}=0$. Knowing that $\epsilon_i \neq 0\;\;i\geq 1$, it comes down to the study of the sequence
$$S_{i,j,k}=\underset{t=\infty}{\hbox{Res}}\; (\epsilon_i^{-1} Q_i(t)+\alpha P_i)(\epsilon_j^{-1}Q_j(t)+\alpha P_j)(\epsilon_k^{-1}Q_k(t)+\alpha P_k)(t^2-1)^2 $$
We have moreover that the system $Sys$ vanish for $P_n$ and $\epsilon_n^{-1} Q_n$, and so is also vanishing for $\epsilon_n^{-1}Q_n +\alpha P_n$. Thanks to that, we will be able to find a recurrence on $S_{i,j,k}$, and moreover it will not depend on $\alpha$.

\begin{lem}
The sequence $S_{i,j,k}$ satisfy the following recurrence relation
\begin{equation}\label{eq4}\begin{split}
-(1+i)(i+j+k+2)S_{i,j,k}+4i(i-1)(i-2)(i-3-j-k)S_{i-2,j,k}+\\
4i(2i-1)j(j-1)S_{i-1,j-1,k}+4i(2i-1)k(k-1)S_{i-1,j,k-1}=0
\end{split}\end{equation}
We have also that $S_{i,j,k}$ is invariant by permutation of the index, and this recurrence relation is not, and so this produce other recurrence relations.
\end{lem}

This recurrence relation can be proved automatically using the Mgfun package for Maple, the holonomic package for Mathematica, or even at hand using integration by parts and a formula between the derivative of $P_n$ and $P_n,P_{n-1}$. Now we will in fact prove much more than necessary for proving Theorem~\ref{thm6}, we will find closed form solutions for the residue we need to compute.

\begin{lem}
For $k\geq j \geq i\geq 1$, we pose
$$f(i,j,k)=\frac{{2}^{d}i!j!k!\Gamma\left(\frac{1}{2}(d+1)\right)\Gamma\left(\frac{a}{2}\right)\Gamma\left(\frac{b}{2}\right)\Gamma\left(\frac{c}{2}\right)} {\Gamma\left(\frac{1}{2}(a+3)\right)\Gamma\left( \frac{1}{2}(b+3)\right)\Gamma\left(\frac{1}{2}(c+3)\right)\Gamma\left(\frac{1}{2}(d+4)\right)}$$
with $a=-i+j+k$, $b=i-j+k$, $c=i+j-k$, $d=i+j+k$. The expression of $S_{i,j,k}$ is given by
\begin{align*}
S_{i,j,k}=\hbox{eval}\left( S_{i,j,k},\alpha=0\right) +\\
\left\lbrace \begin{array}{c}
\lim\limits_{\epsilon\rightarrow 0} \frac{3}{8\pi}f(i+\epsilon,j+\epsilon,k+\epsilon)\alpha^2 \quad\hbox{ if } i+j+k \hbox{ mod }2=1  \\
\lim\limits_{\epsilon\rightarrow 0} \frac{\pi}{16}\frac{1}{\Gamma(\epsilon)} f(i+\epsilon,j+\epsilon,k+\epsilon)\alpha \quad\hbox{ if } i+j+k \hbox{ mod }2=0  \\
\end{array} \right.
\end{align*}
\end{lem}

The limit can be easily computed for all $i,j,k\in\mathbb{N}^*$ but there are no closed form expression for the limit valid for all $i,j,k$. The limit depend in fact of the order of $i,j,k$. We choose in fact such a complicated formula because of its generality. It holds in all cases and thus allows to speed up the proof, avoiding to do $3$ times the same thing, and show effectively the symmetry between the index. With these formulas, it will be easy to prove Theorem~\ref{thm6} because the cases $\partial_{\alpha} S_{i,j,k}=0$ correspond to singular values of the $\Gamma$ functions in the denominator.

\begin{pf} 

\subsubsection{Case $i+j+k\hbox{ mod } 2=1$}
We begin by looking at $f$ for $i+j+k \hbox{ mod }2=1$. This is the easy case, because when we replace $\epsilon$ by $0$ in $f(i+\epsilon,j+\epsilon,k+\epsilon)$, the expression $f(i,j,k) \;\;i,j,k\in\mathbb{N}^*$ is still meaningfull if we assume $\Gamma(-n)=\infty,\;n\in\mathbb{N}$. Indeed, there can be at most one term of this kind and always in the denominator. The corresponding value of $S_{i,j,k}$ will be $0$. We then check than this expression formally satisfy recurrence~\eqref{eq4}. First of all, we remark that when we select the coefficient $\alpha^2$ in $S_{i,j,k}$, we find
$$\hbox{coeff}(S_{i,j,k},\alpha^2)=3\underset{t=\infty}{Res}\; P_{i}P_{j}P_{k} \arctanh\left(\frac{1}{t}\right)(t^2-1)^2$$
Then we make a series expansion of $\arctanh\left(\frac{1}{t}\right)$ for $t=\infty$, and we also get that (noting that $P_iP_jP_k$ is an even polynomial)
$$\hbox{coeff}(S_{i,j,k},\alpha^2)=\sum\limits_{p=0}^{\infty} \frac{3\;\hbox{coeff}(P_{i}P_{j}P_{k}(t^2-1)^2,t^{2p})}{2p+1} =\frac{3}{2}\int\limits_{-1}^1 P_{i}P_{j}P_{k}(t^2-1)^2 dt$$
This relation is interesting because we have an orthogonality property on the polynomials $P_i$. We now need to study boundary cases. Using the recurrence, we come down to $i=1,2$. The recurrence relation~\eqref{eq4} get simpler for $i=2$ and so we can express $S_{2,j,k}$ in function of $S_{1,j,k}$. Using orthogonality, we get the formulas
$$\int\limits_{-1}^1 P_{j}P_{k}(t^2-1)^2 dt=0 \;\;\quad j-k\neq -2,0,2$$
which give $S_{1,j,k}=0,\;\;j-k\neq -2,0,2$. We get for $f(1,j,k)=$
\begin{align*}
-\frac {48(-1)^{\frac{1}{2}(j-k)}\sin \left( \frac{1}{2}\pi\left( 4+j-k \right)  \right) {2}^{j+k}\Gamma  \left( k+1 \right) \Gamma  \left( j+1 \right) }{ \left(2+j-k \right)  \left( j-k \right)  \left( j-k-2 \right) \pi \, \left( 3+j+k \right) \left( j+1+k \right)  \left(j+k-1 \right) }
\end{align*}
Because of the sinus term, this formula vanish for $j-k\neq-2,0,2$. One just need to check it is right in these three cases. We have
$$\int\limits_{-1}^1 P_k^2(t^2-1)^2 dt=\frac{3}{4}{\frac {{4}^{k}\Gamma  \left( k-1/2 \right)\Gamma(k+1)^2}{\Gamma\left(k+5/2\right)}}$$
$$\int\limits_{-1}^1 P_kP_{k+2}(t^2-1)^2 dt=-\frac{3}{2}{\frac {{4}^{k}\Gamma  \left( k+1/2 \right) \Gamma  \left( 3+k \right) \Gamma  \left( k+1 \right) }{\Gamma \left( k+7/2 \right) }}$$
using Mgfun, and we check these are the same as $f$.

\subsubsection{Case $i+j+k\hbox{ mod } 2=0$}
We now look at the function $f$ for $i+j+k \hbox{ mod }2=0$. This time, if we replace formally $\epsilon=0$ in $\frac{1}{\Gamma(\epsilon)} f(i+\epsilon,j+\epsilon,k+\epsilon)\alpha$, we find a quotient of $\Gamma$ functions and in the numerator at most a term of the form $\Gamma(-n),\;n\in\mathbb{N}$. We can still regularize the formula using the relation $\Gamma(n+1)=n\Gamma(n)$. We get in particular that the limit
$$\lim\limits_{\epsilon\rightarrow 0} \frac{1}{\Gamma(\epsilon)} f(i+\epsilon,j+\epsilon,k+\epsilon)$$
is always finite. If there are no term of the form $\Gamma(-n),\;n\in\mathbb{N}$ in the numerator of $f$, then the limit is zero. Using invariance by permutation, we can suppose that $k\geq j\geq i$, and so the only possible infinite term in the numerator is the term in $i+j-k$. We get then a zero limit for $k<i+j$, and for $k\geq i+j$, we can regularize the formula.

We then check that the formula satisfy the recurrence. Now let us look at the boundary cases. With the recurrence, we can reduce $i$ and then $S_{i,j,k}$ is completely determined by $S_{i,j,k},\;i=1,2$. Looking at relation~\eqref{eq4} for $i=2$, it simplifies and we can express $S_{2,j,k}$ in function of $S_{1,j,k}$. Using relation~\eqref{eq4} after permutating the index, we get the recurrence
$$-(1+i)(i+j+3)S_{1,i,j}+4i(i-1)(i-2)(i-j-4)S_{1,i-2,j}+$$
$$4i(2i-1)j(j-1)S_{1,i-1,j-1}=0$$
which allows us to reduce $j$ to $1$. So we only need to compute $S_{1,1,k}$ using Mgfun
$$S_{1,1,2+2k}=-8\,{\frac {16^k\Gamma(k+3/2)\Gamma(k-1/2)\Gamma(k+2)\Gamma(k+1)}{\Gamma  \left( k+4 \right) \Gamma  \left( k+5/2 \right) \sqrt {\pi }}}$$
Our expression correspond to this one.
\end{pf}

\bigskip

To conclude the proof of table $A$ for non zero index, we now look to the expression of $f$. Looking at the formula of $S$ for $i+j+k \hbox{ mod }2=1$, we see that it vanishes exactly when one of the three quantities
$$-i+j+k+3\;\;\;-i+j+k+3\;\;\;-i+j+k+3$$
is non positive. This exactly correspond to the formulas of table $A$ for $i+j+k\hbox{ mod } 2=1$. In the case $i+j+k\hbox{ mod } 2=0$, the quantity
$$\lim\limits_{\epsilon\rightarrow 0} \frac{1}{\Gamma(\epsilon)} f(i+\epsilon,j+\epsilon,k+\epsilon)$$
vanish if and only if all the numbers $-i+j+k,-i+j+k,-i+j+k$ are positiv, which is equivalent using the parity condition to
$$-i+j+k\geq 2\;\;\;-i+j+k\geq 2\;\;\;-i+j+k \geq 2$$

\subsubsection{Case of a zero index}

We now look to the case with at least one zero index. We can write
$$P_0(t)=\frac{t}{t^2-1}\quad Q_0(t)=\frac{1}{t^2-1}\quad Q_i(t)=\epsilon_i P_i(t) \arctanh\left(\frac{1}{t}\right)+\frac{W_i(t)}{t^2-1}$$
(the notation $P_0$ and $Q_0$ is arbitrary here because they are both rational). We begin by the case where exactly one index is zero. We need to compute the residues
$$\underset{t=\infty}{\hbox{Res}}\; (\epsilon_i^{-1}Q_i(t)+\alpha P_i)(\epsilon_j^{-1}Q_j(t)+\alpha P_j)(t^2-1)$$
$$\underset{t=\infty}{\hbox{Res}}\; (\epsilon_i^{-1}Q_i(t)+\alpha P_i)(\epsilon_j^{-1}Q_j(t)+\alpha P_j)t(t^2-1)$$
These are polynomials of degree at most $2$ in $\alpha$ but the coefficient in $\alpha^2$ is always zero because we take the residue of a polynomial at infinity. So one just need to compute the residue in $\alpha$. We expand and suppress the polynomial terms and we get the formula
$$\underset{t=\infty}{\hbox{Res}}\; (\epsilon_i \epsilon_j^{-1}+\epsilon_j\epsilon_i^{-1})\arctanh\left(\frac{1}{t}\right) P_i P_j(t^2-1)$$
Knowing that $\epsilon_i \epsilon_j^{-1}+\epsilon_j\epsilon_i^{-1}$ does not vanish, we can suppress it and we only need to compute the following sequences
\begin{align*}
S^1_{i,j}= \underset{t=\infty}{\hbox{Res}}\;\arctanh\left(\frac{1}{t}\right) P_i P_{i+j}(t^2-1)\\
S^2_{i,j}= \underset{t=\infty}{\hbox{Res}}\;\arctanh\left(\frac{1}{t}\right) P_i P_{i+j} t(t^2-1)
\end{align*}
One just need to prove that either $S_{i,j}^1$ or $S_{i,j}^2$ is not zero for $i\in\mathbb{N}^*$, $j=0,1$ (the condition on the index of table $A$ such that $A=0$ correspond here to $-2\leq i-j\leq 2$ and we use symmetry of the index). Using the parity on $t$ of the polynomials $P_i$, we find that only $S^1_{i,0},S^2_{i,1}$ can be non zero. These sequences can be easily computed for finding a recurrence with Mgfun and then a closed form
$$S^1_{i,0}=-2\frac{4^i\Gamma(i+1)^2}{(i+1)(2i+1)i}\qquad S^2_{i,1}=-4\frac{4^i\Gamma(i+1)^2}{(2i+1)(2i+3)}$$
These expressions do not vanish.
\end{pf}

\subsubsection{Integrability in the cases where $A_{i,j,k}=1$}
Second part: We now prove that if all the non zero terms of second order variational equation correspond only to cases such that $A_{i,j,k}=1$, then the Galois group is abelian. We use the following lemma

\begin{lem}\label{thmlem}
We consider $F(t)=$
$$H_3(t)  \arctanh\left(\frac{1}{t}\right)^3+H_2(t)  \arctanh\left(\frac{1}{t}\right)^2+H_1(t)  \arctanh\left(\frac{1}{t}\right)+H_4(t) $$
with $H_1,H_2,H_3,H_4\in\mathbb{C}[t]$. If the conditions of Theorem~\ref{thm2} are satisfied, then
\begin{itemize}
 \item If $\underset{t=\infty}{\hbox{Res}}\; F(t) =0$, then $\int F\;dt \in\mathbb{C}\left[t,\arctanh\left(\frac{1}{t}\right) \right]$
 \item If $\underset{t=\infty}{\hbox{Res}}\; F(t) \neq 0$, then $\int F\;dt \in\mathbb{C}\left[t,\arctanh\left(\frac{1}{t}\right),\ln\left( t^2-1 \right) \right]$
\end{itemize}
\end{lem}

\begin{pf}
We proceed using integration by parts. We derive the term in $\arctanh\left(\frac{1}{t}\right)^3$. Posing
$$J(t)=\int\limits_{-1}^t H_3(s) ds$$
we have that $J(1)=0$ using the condition in $\alpha^2$ of Theorem~\ref{thm2} and making a series expansion at infinity of $\arctanh\left(\frac{1}{t}\right)$. Then $(t^2-1)$ divide the polynomial $J$. After integration by parts, we get a term of the form
$$-3R(t) \arctanh\left(\frac{1}{t}\right)^2$$
with $R$ a polynomial. Let us try another integration by parts. We get the term
$$-\frac{2}{t^2-1}\int -\frac{3J(t)}{t^2-1}+H_2(t) dt \arctanh\left(\frac{1}{t}\right)$$
We want that this term can be written $Z(t)\arctanh\left(\frac{1}{t}\right)$ with $Z$ a polynomial (with a good choice of integration constant). We only need that
$$\int\limits_{-1}^1 -\frac{3J(t)}{t^2-1}+H_2(t) dt=0$$
Let us look now at the coefficient in $\alpha$ of the residue~\eqref{eq2}. We know it is equal to zero.
$$\hbox{coeff}\left(\underset{t=\infty}{\hbox{Res}}\; F\left(t,\arctanh\left(\frac{1}{t}\right)+\alpha \right) ,\alpha\right)=$$
$$\underset{t=\infty}{\hbox{Res}}\; 3H_3(t) \arctanh\left(\frac{1}{t}\right)^2+2H_2(t) \arctanh\left(\frac{1}{t}\right)$$
which gives using an integration by part (we can see the residue as an integration along a small circle around infinity)
$$\underset{t=\infty}{\hbox{Res}}\; -\frac{6 J(t)}{t^2-1} \arctanh\left(\frac{1}{t}\right)+2H_2(t) \arctanh\left(\frac{1}{t}\right) dt$$
Using the Taylor expansion of $\arctanh\left(\frac{1}{t}\right)$ at infinity, we get
$$\frac{1}{2}\int\limits_{-1}^1 -\frac{6J(t)}{t^2-1}  +2H_2(t) dt=0 $$
This is exactly our condition~\eqref{eq2}. So the last remaining integral to compute is of the form
$$\int Z(t)\arctanh\left(\frac{1}{t}\right) dt \in\mathbb{C}\left[t,\arctanh\left(\frac{1}{t}\right),\ln\left( t^2-1 \right) \right]$$
which can be proved using an integration by part. Now let us look closer to the possible terms in $\ln\left( t^2-1 \right)$. Let suppose that is exists a term $\ln\left( t^2-1 \right)$ in $\int F\;dt$. We will have
$$\int F dt=Z_3(t) \arctanh\left(\frac{1}{t}\right)^3+\dots + Z_0(t) + r \ln\left( t^2-1 \right)$$
with $Z_3,\dots,Z_0$ polynomials and $r$ is a constant because $\ln\left( t^2-1 \right)$ do not appear in $F$. A function in $\mathbb{C}\left[t,\arctanh\left(\frac{1}{t}\right) \right]$ is meromorphic near infinity. Derivating this expression will give
$$F=g' +\frac{rt}{t^2-1} $$
with $g$ a meromorphic function on a neighborhood of infinity. Then
$$\underset{t=\infty}{\hbox{Res}}\; F(t) =r$$
So, if this residue is zero, there will be no $\ln\left( t^2-1 \right)$ terms.
\end{pf}

\bigskip

The integrals to compute for the solutions of second order variational method are the following
$$\int (t^2-1)^2Q_i(t)Q_j(t)Q_k(t) dt \qquad \int (t^2-1)^2P_i(t)Q_j(t)Q_k(t) dt $$
$$ \int (t^2-1)^2P_i(t)P_j(t)Q_k(t) dt\qquad \int (t^2-1)^2P_i(t)P_j(t)P_k(t) dt$$
They are all of the form given by Theorem~\ref{thm2}. We already know that the first one, the third one and the last one satisfy the condition of Theorem~\ref{thm2}, so they all belong to
$$\mathbb{C}\left[t,\arctanh\left(\frac{1}{t}\right),\ln\left( t^2-1 \right) \right]$$ 
For the second one, we compute the coefficient in $\alpha$ of the residue. This gives
$$\hbox{coeff}(\underset{t=\infty}{\hbox{Res}}\; (t^2-1)^2P_i(t)Q_j(t)Q_k(t) dt,\alpha) =\frac{1}{2}\epsilon_j\epsilon_k \int\limits_{-1}^1 P_i(t)P_j(t)P_k(t)(t^2-1)^2 dt$$
which equals to zero because it corresponds to the condition in $\alpha^2$ for the first one. So the residue condition is also satisfied, and then thanks to Theorem~\ref{thm2}, it also belongs to
$$\mathbb{C}\left[t,\arctanh\left(\frac{1}{t}\right),\ln\left( t^2-1 \right) \right]$$

\end{pf}

\subsection{Applications}

\subsubsection{The non diagonalizable case}
Let us look first at this diagonalizability condition. We will see it is in fact not a strong condition, because if the potential is integrable at order one, there is only one possibility for which the Hessian matrix can be non diagonalizable, only for the eigenvalue $-1$. We do not make a complete analysis of this case at order $2$ because it is not possible to make an efficient reduction in this case to produce a nice criterion and this case is very rare in practice (still in this case, the Theorem~\ref{thmmain1} gives integrability constraints on a subsystem, but this constraint is not a priori optimal), so it is probably more adapted to make an analysis case by case directly in applications.

\begin{thm}\label{thm8} We consider the equation
\begin{equation}\label{eq5}
\ddot{X(t)}=\frac{1}{\phi(t)^3} AX(t)
\end{equation}
with $A$ a matrix in Jordan form. Then this equation has a virtually abelian Galois group if and only if
\begin{itemize}
\item $A_{i,i}=\frac{(p_i-1)(p_i+2)}{2}$ with $p_i\in\mathbb{N}$
\item $A$ is diagonal except maybe for eigenvalue $-1$ for which the Jordan blocks should have a size less than $2$.
\end{itemize}
\end{thm}

\begin{pf}
In the non diagonalizable case, the non homogeneous part of the variational equation correspond to terms outside the diagonal in the matrix $A$. If the Jordan block have size less tan $2$, the equation~\eqref{eq5} can be rewritten
\begin{equation}\label{eq6}
\ddot{X}=\frac{1}{\phi(t)^3} d_i X +\frac{1}{\phi(t)^3}Y\qquad \ddot{Y}=\frac{1}{\phi(t)^3} d_i Y
\end{equation}
We have in particular that $Y$ satisfy the homogeneous equation for the same eigenvalue, because a Jordan block has the same eigenvalue on the diagonal. For the case with bigger Jordan block, we would get even stronger conditions because then the equation~\eqref{eq5} would be a subsystem of equation~\eqref{eq6}. In our case (except for eigenvalue $-1$), this analysis is not necessary. We compute the solutions and we prove that the following function should be in the Picard Vessiot field.
$$S=\int (t^2-1)(Q_i(t)+\epsilon_i \alpha P_i)^2 dt$$
With Theorem~\ref{thm2} and Lemma~\ref{thmlem}, we know it is enough to study the sequence
$$S_i=\underset{t=\infty}{\hbox{Res}}\; (t^2-1)(Q_i(t)+\epsilon_i \alpha P_i)^2$$
The interesting term is in $\alpha$, because the coefficient in $\alpha^2$ is always zero. The Mgfun package give us a recurrence and then a closed form for this residue
$$S_i=-2\frac{4^i i \Gamma(i)^2}{(2i+1)(i+1)}$$
which is never zero, except for $i=0$. The case $i=0$ necessitate the analysis of a Jordan block of size $3$ (higher Jordan block size would still have this system as a subsystem). This gives the equation
$$\ddot{X_1}=-\frac{1}{\phi(t)^3} X_1\qquad \ddot{X_2}=\frac{1}{\phi(t)^3}(- X_2+X_1) \qquad \ddot{X_3}=\frac{1}{\phi(t)^3}(-X_3 +X_2)$$
Again, we compute the solution and we use Theorem~\ref{thm2} to prove non commutativity of the monodromy (in this case, it is very easy because there are no parameters). To conclude, we notice that the Galois group is always connected because we only take recursively integrations (no algebraic functions are involved) and because the Galois group of $\ddot{X}=\frac{1}{\phi(t)^3} d_i X$ is always connected.
\end{pf}

\subsubsection{A useful corollary}

\begin{thm}\label{thmappli}
Let $V$ be a meromorphic homogeneous potential of degree $-1$ in dimension $n$, $c$ a Darboux point of $V$ with multiplier $-1$. We pose $\lambda_i\;\;i=1\dots n$ the eigenvalues of $\nabla^2V(c)$ with $\lambda_1=2$ (the eigenvalue $2$ always appear in the spectrum). Suppose that $\nabla^2V(c)$ is diagonalizable and
$$\{\lambda_2,\dots,\lambda_n \} \subset \{ (2k-1)(k+1),\; k\in B \}$$
where $B\subset \mathbb{N}$ such that
$$\max(B) \leq \max(2\min(B)-1,0)$$
Then the variational equation at order $2$ near the homothetic orbit associated to $c$ has a virtually abelian Galois group.
\end{thm}

\begin{rem}
In practice, this corollary says us that if the eigenvalues of the Hessian matrix have all an even index and are sufficiently near to each other, then the system is always integrable at order $2$ without any additional conditions. Moreover, this theorem is in some sense "optimal", because if it is not satisfied, then there will be strong additional integrability constraint (of codimension at least $1$) for integrability. It allows also to have a strong intuition about what will be the easy and the hard cases in proving non integrability of a particular problem depending on parameters.
\end{rem}

\begin{pf}
We have the Euler relation due to homogeneity and the Darboux point condition
\begin{equation}\label{eq7}
\sum\limits_{i=1}^n q_i \frac{\partial }{\partial q_i} V =-V \qquad \frac{\partial }{\partial q_i} V(c)=-c_i
\end{equation}
By derivating the Euler relation in $q_j$, we get
$$\sum\limits_{i=1}^n q_i \frac{\partial }{\partial q_i\partial q_j} V =-2\frac{\partial }{\partial q_j}V$$
With the Darboux point relation, this implies that $c$ is an eigenvector with eigenvalue $2$. Let us note $X_1=c,X_2,\dots,X_n$ a basis of eigenvectors of $\nabla^2 V(c)$. We will first prove that
$$D^3(V)(c).(X_1,X_a,X_b)=0\;\;\forall a\neq b$$
We derive the Euler relation \eqref{eq7} two times and evaluate on $c$
$$\sum\limits_{i=1}^n c_i \frac{\partial }{\partial q_i\partial q_j\partial q_k} V(c) =-3\frac{\partial }{\partial q_j\partial q_k}V(c)\qquad \forall j,k$$
We multiply each line with index $j$ by $(X_a)_j$ and we sum over the $j$
$$\sum\limits_{j=1}^n\sum\limits_{i=1}^n (X_a)_j c_i \frac{\partial }{\partial q_i\partial q_j\partial q_k} V(c) =-3\lambda_a (X_a)_k\qquad \forall k$$
with $\lambda_a$ the eigenvalue associated to $X_a$ and using the fact that $X_a$ is an eigenvector of $\nabla^2 V(c)$. We then multiply each line with index $k$ by $(X_b)_k$ and we sum over the $k$
\begin{equation}\label{eq8}
\sum\limits_{k=1}^n\sum\limits_{j=1}^n\sum\limits_{i=1}^n (X_b)_k(X_a)_j c_i \frac{\partial }{\partial q_i\partial q_j\partial q_k} V(c) =-3\lambda_a <X_a | X_b>=0
\end{equation}
thanks to orthogonality. This is the expression of $D^3(V)(c).(X_1,X_a,X_b)$. Let us now use Theorem~\ref{thmmain1}. We first remark that all invoked index of table $A$ are even. Moreover, we have that for three even index $i,j,k$, if $\max(i,j,k) \leq 2\min(i,j,k)-2$, then
$$A_{i,i,i},A_{j,j,j},A_{k,k,k},A_{i,j,j},A_{i,k,k},A_{j,i,i},A_{j,k,k},A_{k,i,i},A_{k,j,j},A_{i,j,k}=1$$
We also have $A_{0,0,0}=1$. So, if the eigenvalues of $\nabla^2 V(c)$ satisfy
$$\{\lambda_1,\lambda_2,\dots,\lambda_n \} \subset \{ (2k-1)(k+1),\; k\in B \}$$
where $B\subset \mathbb{N}$ and such that
$$\max(B) \leq \max(2\min(B)-1,0)$$
then the system is integrable at order $2$. Still, knowing that the eigenvalue $2$ always appear in the spectrum of $\nabla^2 V(c)$, this would be useless. But we know that $c$ is always an eigenvector with eigenvalue $2$, and we have that all the possible conditions linked to this eigenvector are of the form $D^3(V)(c).(X_1,X_a,X_b)=0$. These are automatically satisfied for $a\neq b$. For $a=b$, we have that $A_{2,i,i}=1 \;\forall i\neq 1$, so the only possible problem would be if $X_a$ has the eigenvalue $0$, but in this particular case, we also get with equation \eqref{eq8}
$$\sum\limits_{k=1}^n\sum\limits_{j=1}^n\sum\limits_{i=1}^n (X_a)_k(X_a)_j c_i \frac{\partial }{\partial q_i\partial q_j\partial q_k} V(c) =-3\lambda_a <X_a | X_a>=0$$
because $\lambda_a=0$. So the possible integrability conditions involving the eigenvector $X_1=c$ are always satisfied, and thus we can remove one time the eigenvalue $2$ from $B$. This gives the theorem.
\end{pf}

\subsection{Study of the Galois group in the integrable case}

We now prove Theorem~\ref{thmmain2}, analyzing more precisely the Galois group in the case where the integrability conditions of Theorem~\ref{thm2} are satisfied. We will see that in fact the Galois group almost never grows, and the Galois group can be in fact precisely computed thanks to Lemma~\ref{thmlem}.

\begin{pf}
In the integrable case, the variational equations at order $1$ involve the functions $P_i,Q_i$ which are in $\mathbb{C}(t,\arctanh\left(\frac{1}{t}\right))$ (after variable changes). The only non multivalued function $Q$ is the function $Q_0$, but the eigenvalue $2$ is always in the spectrum, and so the Galois group is always $\mathbb{C}$.\\
At order $2$, using Theorem~\ref{thm2}, we already know that the solutions are in $\mathbb{C}\left[t,\arctanh\left(\frac{1}{t}\right),\ln\left( t^2-1 \right) \right]$ and we know a condition for which the term in $\ln\left( t^2-1 \right)$ do not appear. Thanks to Lemma~\ref{thmlem}, we know that this logarithmic term can appear only if $S_{i,j,k}$ is a non zero constant (independent of $\alpha$ because we suppose that the second order variational equations have a virtually abelian Galois group). Let us prove that
\begin{equation}\label{eq9}
\hbox{eval}(S_{i,j,k},\alpha=0)=0\;\;\forall \;i,j,k \in\mathbb{N}^*
\end{equation}
We only need to use the recurrence \eqref{eq4} for $S_{i,j,k}$. To prove that this sequence is zero, we then only need to prove it vanishes on the boundary, and here it comes down to the case $i=j=1$ (because of the invariance by permutation and the recurrence on the plane $i=1$). Mgfun then prove that
$$\hbox{eval}(S_{1,1,k},\alpha=0)=0\;\;\forall \;k \in\mathbb{N}^*$$
With recurrence \eqref{eq4}, we get the property \eqref{eq9}.

\medskip

Let us look now to the case where one of the index is zero. We need to study
$$S^1_{i,j}=\underset{t=\infty}{\hbox{Res}}\; (t^2-1)Q_i(t)Q_j(t)\qquad S^2_{i,j}=\underset{t=\infty}{\hbox{Res}}\; (t^2-1)tQ_i(t)Q_j(t)$$
We also prove they vanish for $\alpha=0$ using recurrence. If two index are zero, then we need to study
$$S^1_{i}=\underset{t=\infty}{\hbox{Res}}\; Q_i(t)\qquad S^2_{i}=\underset{t=\infty}{\hbox{Res}}\; tQ_i(t) \qquad S^3_{i}=\underset{t=\infty}{\hbox{Res}}\; t^2Q_i(t)$$
All these sequences are zero except for $i=1$ for which $S^3_1=-\frac{2}{3}$. Eventually, in the case where all index are zero, we need to compute the following integrals
$$\int \frac{1}{t^2-1} dt\qquad \int \frac{t}{t^2-1} dt\qquad \int \frac{t^2}{t^2-1} dt\qquad \int \frac{t^3}{t^2-1} dt$$
The second and the forth integral have a term in $\ln\left( t^2-1 \right)$. This gives Theorem~\ref{thmmain2}.
\end{pf}

\begin{rem}
The computation of the sequence $\hbox{eval}\left(S_{i,j,k},\alpha=0\right)$ has in fact no sense when $S_{i,j,k}$ depend on $\alpha$. Indeed $\alpha$ corresponds to the multivaluation $Q_i$. If $\hbox{eval}\left(S_{i,j,k},\alpha=0\right)=0$ and $S_{i,j,k}$ depend on $\alpha$, then when we replace $\alpha$ by $\alpha+1$, we get $\hbox{eval}\left(S_{i,j,k},\alpha=0\right)\neq 0$ and so this vanishing term correspond in fact only to a convention taken for the $Q_i$. Still the convention is well chosen, because it allows to study $\hbox{eval}\left(S_{i,j,k},\alpha=0\right)$ without making a distinction between integrable cases and non integrable cases (in particular, $\hbox{eval}\left(S_{i,j,k},\alpha=0\right)$ is almost always zero for all values, and this property allows a much faster proof than in previous sections).
\end{rem}

\label{}\bibliographystyle{model1-num-names}
\bibliography{nonzeroangular}

\end{document}